\documentclass{article}
\usepackage[leqno]{amsmath}
\usepackage{amsthm,amscd}
\usepackage{diagrams}
\usepackage{amsfonts} 
\usepackage{amssymb} 

\newtheorem{theorem}{Theorem}[section]
\newtheorem{definition}[theorem]{Definition}
\newtheorem{example}[theorem]{Example}
\newtheorem{lemma}[theorem]{Lemma}
\newtheorem{proposition}[theorem]{Proposition}
\newtheorem{corollary}[theorem]{Corollary}

\newcommand{\A}{{\mathcal A}}
\newcommand{\B}{{\mathcal B}}

\newcommand{\R}{{\mathbb R}}

\newcommand{\bbS}{{\mathbb S}}

\newcommand{\codim}{\operatorname{codim}}
\newcommand{\im}{\operatorname{im}}

\newcommand{\ch}{{\rm \bf Ch}}

\newcommand{\chhs}{\epsilon_{j}^{\sigma}}
\newcommand{\phhs}{\varphi_{j}^{\sigma}}
\newcommand{\phis}{\varphi_{j}^{\sigma}}  
\usepackage{graphicx}

\begin{document}
%
%

\title{Chambers of Arrangements of  Hyperplanes and Arrow's Impossibility Theorem}

\author{Hiroaki Terao
\thanks{This research was
supported in part by Japan Society for the Promotion
of Science.}
\\
{\it Hokkaido University, Sapporo, Japan} 
\\
}

\date{\today}

\maketitle

\begin{abstract}
Let ${\mathcal A}$ be a nonempty real central arrangement of hyperplanes
and ${\rm \bf Ch}$ be the set of chambers of ${\mathcal A}$.  
Each hyperplane
$H$ defines a half-space $H^{+} $ and 
the other half-space
$H^{-}$.
Let $B = \{+, -\}$. 
For $H\in {\mathcal A}$, 
define a map
$\epsilon_{H}^{+}  : {\rm \bf Ch} \to B$ 
by
$
\epsilon_{H}^{+}  (C)
=
+ ~\text{(if~} C\subseteq H^{+})
\,
\text{~and~} 
\epsilon_{H}^{+}  (C)
= -  ~\text{(if~} C
\subseteq H^{-}).
$
Define
$
\epsilon_{H}^{-}
=
-\epsilon_{H}^{+}.
$
Let ${\rm \bf Ch}^{m} = {\rm \bf Ch}
\times{\rm \bf Ch}\times\dots\times{\rm \bf Ch}
\,\,\,
(m\text{~times}).$  
Then 
the maps 
$\epsilon_{H}^{\pm}$ induce
the maps
$\epsilon_{H}^{\pm} : {\rm \bf Ch}^{m}  \to B^{m} $. 
We will study the admissible maps
$\Phi
:
{\rm \bf Ch}^{m} \to {\rm \bf Ch}$ which are compatible with every
$\epsilon_{H}^{\pm}$.
Suppose 
$|{\mathcal A}|\geq 3$
and
$m\geq 2$. 
Then we will show that 
${\mathcal A}$ is indecomposable
if and only if 
every admissible map
is a projection to a component.
When ${\mathcal A}$ is a braid arrangement, which is 
indecomposable, this 
result is equivalent to Arrow's
impossibility theorem in economics.
We also determine the set of admissible maps
explicitly
for every nonempty
real central arrangement.

{\bf Key words:}
arrangement of hyperplanes, 
chambers,
braid arrangements,
Arrow's impossibility theorem.
\end{abstract}

\section{Main Results}
Let $\A
=
\{
H_{1}, H_{2}, \dots, H_{n}   
\}
$  be a nonempty
real central arrangement of hyperplanes
in 
$\R^{\ell}$.
In other words,
each hyperplane $H_{j} $ goes through the origin of
$\R^{\ell}$.
In this note,
we frequently refer to
\cite{OrT}
for elementary facts about arrangements of hyperplanes,
which are usually referred as 
{\bf arrangements} for brevity.
The connected components of the complement
$ \R^{\ell} \setminus \bigcup_{1\leq j\leq n} H_{j} $
are called {\bf chambers}
of $\A$.
Let $\ch = \ch(\A)$ denote the set of chambers
of $\A$.
For each hyperplane $H_{j} \in\A$, 
fix a
real linear form $\alpha_{j} $ such that
$H_{j} =\ker(\alpha_{j} ).$
The product 
$\prod_{j=1}^{n} \alpha_{j}  $ is called a 
{\bf defining polynomial}
for $\A $.
Define
\[
H_{j}^{+} = \{x\in\R^{\ell} \mid \alpha_{j}(x) > 0\},
\,\,\,\,
H_{j}^{-} = \{x\in\R^{\ell} \mid \alpha_{j}(x) < 0\}
\,\,\,
(j = 1,\dots, n).
\]
Throughout this note, let $\sigma$ denote $+$ or $-$.
Let $B=\{+, -\}$, which
we frequently 
consider as a multiplicative group of order
two in the natural way.

Let $1\leq j\leq n$. 
The maps
$\chhs : \ch \longrightarrow B$
are defined by
$
\chhs(C) = \sigma\tau
$
if $C\subseteq H_{j}^{\tau}$
$(\sigma, \tau \in B)$.  
Let $m$ be a positive integer.
Consider the $m$-time direct products
$\ch^{m}$ and $B^{m} $. 
We let the same symbol
$\chhs$ also denote the map
$\ch^{m} \to B^{m} $ induced from
$\chhs
:
\ch \to B$:

\[
\chhs
(C_{1} , C_{2} , \dots , C_{m} ) 
=
(
\chhs(C_{1}) , 
\chhs(C_{2}) , 
\dots , 
\chhs(C_{m}) 
)
\]
for $
(C_{1} , C_{2} , \dots , C_{m} ) 
\in\ch^{m}.$

\begin{definition}
\label{admissiblemaps}
A map 
$\Phi :
\ch^{m} \longrightarrow \ch$
is called an {\bf admissible map} 
if 
 there exists a family of maps
$\phhs
:
B^{m} \longrightarrow B$
$(1\leq j\leq n, \,\,\,
\sigma\in B=\{+, -\})$
which satisfies the following two conditions:

$(1)$ 
$\phhs (+,+,\dots, +) = +$,
and

$(2)$
the diagram
\begin{diagram}
\ch^{m}& \rTo^{\Phi}  &\ch\\
\dTo^{\chhs} &  &\dTo_{\chhs}\\
B^{m}& \rTo{\phhs}  &B
\end{diagram}
commutes
for each $j, 1\leq j\leq n$, and 
$\sigma\in B=\{+, -\}$. 

Let
$AM(\A, m)$
denote
the set of all
admissible maps determined by 
$\A$ and $m$.
\end{definition}

\smallskip

As we will see in Proposition
\ref{uniquely}, 
when $\Phi$ is an admissible map,
a family of 
maps
$\phhs$ $(1\leq j\leq n, \sigma\in B= \{+, -\})$ 
satisfying the conditions in Definition \ref{admissiblemaps} 
is uniquely determined by $\Phi$, $\A$ and $m$.

The main purpose of this note is to
study the set $AM(\A, m)$ for all
$\A$ and $m$.

\begin{definition}
For $ 1\leq h\leq m$,
let
\begin{align*} 
\Phi
&= \text{the projection to the $h$-th component},\\
\phhs
&=\text{the projection to the $h$-th component}.
\end{align*} 
Then it is easy to see that 
$\Phi$ is an admissible map
with a family of maps $\phhs\,\,(1\leq j\leq n,
\sigma\in B=\{+, -\})$.
We call the admissible maps of this type
{\bf projective admissible maps.}
\end{definition}

For a central arrangement $\A$, define
$$
r(\A) 
=
\codim_{\R^{\ell} }  \bigcap_{1\leq j\leq n} H_{j}.$$  

\begin{definition}
A central arrangement 
$\A$ is said to be 
{\bf decomposable} if there exist
nonempty arrangements $\A_{1} $ and
$\A_{2} $ such that
$\A = \A_{1} \cup \A_{2} $ (disjoint)
and
$r(\A) = r(\A_{1} ) + r(\A_{2} )$. 
In this case, write
$\A = \A_{1} \uplus \A_{2} $.
A central arrangement 
$\A$ is said to be 
{\bf indecomposable} if it is not decomposable.
\end{definition}

Note that 
$\A = \A_{1} \uplus \A_{2} $
if and only if
the
defining polynomials for $\A_1$ and $\A_2$ have 
no common variables
after an appropriate 
linear coordinate change.  

\medskip

\noindent
{\bf Remark.}
It is also known \cite[Theorem 2.4 (2)]{STV}
that $\A$ is decomposable if and only
if
its Poincar\'e polynomial 
\cite[Definition 2.48]{OrT}
$\pi(\A, t)$ is divisible by
$(1+t)^{2}$.  

\medskip

We will see in Proposition \ref{decomp} that
any nonempty
real central arrangement 
$\A$ can be uniquely (up to order)
decomposed into nonempty
indecomposable arrangements$:$
$$
\label{bunkai} 
\A = \A_{1} \uplus \A_{2} \uplus \dots \uplus \A_{r}. 
\eqno{(*)}
$$

The following two theorems
completely determine 
the set
$AM(\A, m)$ of admissible maps.

\begin{theorem}
\label{maindecomp} 
For a nonempty
real central arrangement 
$\A$ with the 
decomposition $(*),$
there exists a natural bijection
\[
AM(\A, m) \simeq AM(\A_{1}, m) \times AM(\A_{2}, m) \times \dots 
\times
AM(\A_{r}, m)
\]
for each positive integer $m$. 
\end{theorem}

\begin{theorem}
\label{main} 
Let $\A$ be a nonempty indecomposable
real central arrangement
and $m$ be a positive integer.
Then,

$(1)$ 
if $|\A|=1$,
\[
AM(\A, m) =
\{
\Phi : \ch^{m} \rightarrow \ch \mid 
\Phi(C, C, \dots , C) = C
\text{~for each chamber~ } C\},
\]
 
$(2)$ 
if $|\A|\geq 3$,
every admissible map is projective.

$($Note that, if $|\A| = 2$, then $\A$ is decomposable.$)$ 
\end{theorem}

\begin{corollary}
\label{corollary1}
Decompose
a
nonempty
real central arrangement 
$\A$ into nonempty
indecomposable arrangements
as
\[
\A = \A_{1} \uplus \A_{2} \uplus \dots \uplus \A_{a}
\uplus
\B_{1} \uplus \B_{2} \uplus \dots \uplus \B_{b} 
\]
with
$|\A_{p} | = 1 \,\,(1\leq p \leq a)$ and
$|\B_{q} | \geq 3 \,\,(1\leq q \leq b)$.
Then, for each positive integer $m$, 
\[
|AM(\A, m)| = (2^{a(2^{m}-2)}) m^{b}.  
\]
\end{corollary}

\noindent
{\bf Remark.}
Theorem \ref{main} can be regarded as a generalization
of
Kenneth Arrow's impossibility theorem
(\cite{A, MCWG}) in economics:

In the impossibility theorem,
we assume that a society of $m$ people
have $\ell$ policy options 
and that every individual has his/her own
order of preferences on the $\ell$ 
policy options.
A social welfare function can be interpreted as a voting system
by which the individual preferences are
aggregated into a single societal preference.
We require the following two requirements
for a reasonable social welfare function:

(A) the society prefers the option $i$ 
to the option $j$ if every individual
prefers
 the option $i$ 
to the option $j$
(Pareto property),
and
(B) whether the society prefers
 the option $i$ 
to the option $j$ only depends on which individuals
prefer  the option $i$ 
to the option $j$
(pairwise independence).

The conclusion of Arrow's impossibility theorem
is striking:
for $\ell\geq 3$,
the only social welfare function
satisfying the two requirements (A) and (B)
is a dictatorship, that is, the societal preference has to be
 equal to the
preference of one particular individual. 

In Theorem \ref{main}, let $\A$ 
be a braid arrangement in $\R^{\ell} $
$(\ell \geq 3)$,
i. e., 
$$\A = \{H_{ij} \mid 1\leq i<j\leq \ell\},
\,\,\,
\text{where} 
\,\,\,
H_{ij} := \ker (x_{i} - x_{j}).
$$
The braid arrangements are indecomposable as we will see
in Example
\ref{braid}.
Let 
$H_{ij}^{+}
=
\{(x_{1} , x_{2} , \dots , x_{\ell})\in\R^{\ell} \mid x_{i} > x_{j} \}  $
and
$H_{ij}^{-}
=
\{(x_{1} , x_{2} , \dots , x_{\ell})\in\R^{\ell} \mid x_{i} < x_{j} \}$.
Then each chamber of $\A$ can be uniquely expressed as 
$$\{(x_{1} , x_{2} , \dots , x_{\ell})\in\R^{\ell} \mid 
x_{\pi(1)} 
< 
x_{\pi(2)} 
< 
\dots
<
x_{\pi(\ell)} 
\}
$$
for a permutation $\pi$ of 
$\{1, 2, \dots , \ell\}$.
This gives a one-to-one correspondence between
$\ch(\A)$ and the permutation group $\bbS_{\ell} $  of
$\{1, 2, \dots , \ell\}$.
Thus we can interpret
an 
order of preferences on $\ell$ policy options
as a chamber of a braid arrangement.
Similarly, 
we interpret a social welfare function
as
the map $\Phi$
and the dictatorship by the $h$-th individual
as
the projection to the $h$-th component. 
The 
requirements 
(A) 
(Pareto property)
and (B)
(pairwise independence)
correspond to
the 
conditions 
(1)
($\phhs(+, \dots, +) = +$) 
and 
(2) 
 (commutativity) 
in Definition
\ref{admissiblemaps} 
respectively.
So,
in our terminology, Arrow's impossibility theorem can be formulated
as:

\begin{quotation} 
\noindent
{\it 
If $\A$ is a braid arrangement with $\ell \geq 3$, 
then
every admissible map is projective.}
\end{quotation}

Thanks to Theorems \ref{maindecomp} and \ref{main} 
we have 
the following necessary and sufficient condition
for 
a nonempty real central
arrangement to have the property that
every admissible map is
projective:

\begin{corollary}
\label{corollary2}
Let $\A$ be
a nonempty 
real central arrangement and $m$ be
a positive
integer.
Every admissible map is projective
if and only  
if

$($case 1$)$ $m=1$,
 or

$($case 2$)$ $\A$ is indecomposable with $|\A|\geq 3$.
\end{corollary}

\section{Proof of Theorem \ref{maindecomp} } 
Let $\A
=
\{
H_{1}, H_{2}, \dots, H_{n}   
\}
$  be a nonempty
real central arrangement 
in $\R^{\ell}$.
Let
$\B$ be a subarrangement of $\A$, in other words,
$\B \subseteq \A$.
We say that $\B$  is  {\bf dependent}
if 
$$
r(\B)
=
\codim_{\R^{\ell} } (\bigcap_{H\in \B} H ) < |\B|. $$
A subarrangement $\B$ of $\A$ 
is called {\bf independent}
if it is not dependent.
If $\B$ is a minimally dependent subset, then
$\B$ is called a {\bf circuit.}
If $\B$ is a maximally independent subset
in $\A$, then $\B$ 
is called a {\bf basis}
for $\A$. 

We introduce a graph $\Gamma(\A)$ associated with $\A$.
The set of vertices of  $\Gamma(\A)$ is $\A$.
Two vertices
$H_{j_{1} } , H_{j_{2} } \in\A$ $(j_{1} \neq j_{2} )$ 
 are connected by an edge
if and only if
there exists a circuit (in $\A$)
containing $\{H_{j_{1} } , H_{j_{2} } \}$. 

\begin{lemma}
\label{lemma2.1} 
A nonempty real
central arrangement $\A$ is indecomposable if and only if
the graph $\Gamma(\A)$ is connected.
\end{lemma}

\begin{proof}
If
$\Gamma(\A)$ is disconnected, then decompose $\A$ as
$\A 
=
\A_{1} 
\cup
\A_{2} 
$ 
so that 
$\A_{1} \neq \emptyset$,
$\A_{2} \neq \emptyset$,
 and $\{H_{j_{1} } , H_{j_{2} } \}$ 
is not contained in any circuit
whenever $H_{j_{p} } \in \A_{p} $ 
\,
$
(p = 1, 2).
$
Choose a basis $\B_{p} $ of $\A_{p} $ 
$
(p = 1, 2).
$
Then $\B_{1} \cup \B_{2} $
is also independent
because it does not contain any circuit.
Thus
\[
r(\A) = 
|\B_{1} \cup \B_{2}|
=
|\B_{1}| + |\B_{2}|
=r(\A_{1} )
+
r(\A_{2} ),
\]
 which implies 
$
\A=\A_{1} \uplus \A_{2}. $
So $\A$ is decomposable.

Conversely assume that  
$\A=\A_{1} \uplus \A_{2} $
with 
$\A_{1} \neq \emptyset,$ 
$\A_{2} \neq \emptyset$. 
We may assume,
after an appropriate 
linear coordinate change, that
the
defining polynomials for $\A_1$ and $\A_2$ have 
no common variables.
Let $H_{j_{p}  } \in \A_{p} $ 
$
(p = 1, 2).
$
Suppose that 
there exists a circuit $\B$ containing
$H_{j_{1} } $ and $H_{j_{2} } $.
Then
 $\B \cap \A_{1} $ 
and
 $\B \cap \A_{2} $ 
 are both independent.
This implies that
$\B$ is also independent, which is 
a contradiction.
\end{proof}

\begin{example}
\label{braid}
Let $\A$ be a braid arrangement in $\R^{\ell}
(\ell\geq 2):$
\[
\A = \{
H_{ij} \mid
1\leq i<j\leq \ell
\},
\]
where 
$
H_{ij}
=
\ker(x_{i}-x_{j}).
 $ 
 If $\ell=2$, then $|\A|=1$ and $\A$ 
is indecomposable.  Let $\ell \geq 3$.
Then 
$
\{
H_{ij}, H_{jk}, H_{ik}   
\}
$   
for $1\leq i<j<k\leq \ell$ 
is a circuit.   Thus it is easy to check that
$\A$ is indecomposable by 
applying
Lemma \ref{lemma2.1}. 
\end{example}

By Lemma \ref{lemma2.1}, we immediately
have

\begin{proposition}
\label{decomp}
Any 
nonempty
real central arrangement 
$\A$ can be uniquely
(up to order) decomposed
into nonempty
indecomposable arrangements
\[
\A = \A_{1} \uplus \A_{2} \uplus \dots \uplus \A_{r}. 
\]

\end{proposition}

Let $m$ be a positive integer.
For $S \subseteq \{1, \dots, m\}$,
define
$S_{+} =
 (\sigma_{1}, \dots , \sigma_{m})
\in
B^{m} $
with
\[
\sigma_{i} =
\begin{cases}
+  \,\,\text{if~} i\in S,\\
-  \,\,\text{if~} i\not\in S.
\end{cases}  
\]
Then $(S^{c})_{+}
=
(-\sigma_{1}, \dots , -\sigma_{m})
=-S_{+}$.

\begin{proposition}
\label{surjective}
Assume
$\sigma\in B= \{+, -\}$  and 
$1\leq j\leq n$.
Then the map
$\epsilon_{j}^{\sigma} : \ch^{m} \rightarrow B^{m}  $ is surjective.
\end{proposition}

\begin{proof}
An arbitrary element of $B^{m} $ can be expressed as
$S_{+} $
for
some $S \subseteq \{1,2,\dots, m\}$.
Suppose that
$C$ and $ C'$ 
are chambers such that
$C\subseteq H_{j}^{+}  $ 
and
$C'\subseteq H_{j}^{-}  $.
Define
$
\mathcal C
=
(
C_{1} , C_{2} , \dots , C_{m} 
)\in\ch^{m} 
$  by
\[
C_{i} 
=
\begin{cases}
C \,\,\text{if~} i\in S,\\
C' \,\,\text{if~} i\not\in S.
\end{cases}
\]
Then we have
$
\epsilon_{j}^{+} (\mathcal C) 
=
S_{+}.
$
Let 
$
-\mathcal C
=
(
-C_{1} , -C_{2} , \dots , -C_{m} 
)\in\ch^{m}, 
$ 
where $-C_{i} $ denotes the antipodal chamber of 
$C_{i} $.
Then
$
\epsilon_{j}^{-} (-\mathcal C) 
=
-(S^{c})_{+}
=S_{+}.
$
\end{proof}

\begin{proposition}
\label{uniquely}
When $\Phi$ is an admissible map, a family of maps
$\phhs$ $(1\leq j\leq n, \,\, \sigma\in B = \{+, -\})$ satisfying the
the conditions
in Definition \ref{admissiblemaps} is uniquely determined.
\end{proposition}

\begin{proof}
It is obvious because of
Proposition \ref{surjective}.
\end{proof}

\begin{proposition}
\label{unanimity}
When $\Phi$ is an admissible map, 
$\Phi(C, C, \dots, C) = C$ for any chamber
$C\in\ch$. 
\end{proposition}

\begin{proof}
By Definition \ref{admissiblemaps}, 
two chambers 
$\Phi(C, C, \dots, C)$ and $C$
are on the same side of
every $H_{j} \in\A$.  Thus  
$\Phi(C, C, \dots, C) = C$.
\end{proof}

Suppose that
$\A = \A_{1} \uplus \A_{2}$
with 
$\A_{1} \neq \emptyset$ and
$\A_{2} \neq \emptyset$.   
We may assume that
the
defining polynomials for $\A_1$ and $\A_2$ have 
no common variables.
Then the following lemma is obvious:

\begin{lemma}
\label{Lemma3.1(1)}
The map
\[
\alpha : \ch(\A_{1} )^{m} \times \ch(\A_{2} )^{m} 
\longrightarrow
\ch(\A_{1} \uplus \A_{2} )^{m}, 
\]
given by
$$\alpha(C_{1}, \dots, C_{m}, D_{1}, \dots , D_{m})
=
(C_{1} \cap  D_{1}, \dots, C_{m} \cap D_{m})
$$
for 
$C_{i}\in\ch(\A_{1}), 
D_{i}\in\ch(\A_{2})
\,\,\, 
(
i=1,\dots,m
)
$ , 
is bijective.
\end{lemma}

\smallskip

\begin{lemma}
\label{Lemma3.1(2)}
Let $p \in \{1, 2\}$. 
For $H_{j}\in \A_{p}$, the diagram
\begin{diagram}
\ch(\A_{1} ) \times \ch(\A_{2} )
& 
\rTo^{\alpha}
&
\ch(\A_{1} \uplus \A_{2} )\\
 \dTo^{\pi_{p}} 
&
&
\dTo_{\epsilon_{j}^{\sigma}} \\
\ch(\A_{p})
&
\rTo^{\epsilon_{j, p}^{\sigma}}
&B
\end{diagram} 
is commutative, where 
$\pi_{p} $ is the projection to the $p$-th component, and 
 $\epsilon_{j, p}^{\sigma}$ is the map
$\epsilon_{j}^{\sigma}$ for $\A_{p} $. 
\end{lemma}

\begin{proof}
Let $p=1$ for simplicity.
Then
$$
\epsilon^{\sigma}_{j}\circ \alpha(C, D)
 = 
\epsilon^{\sigma}_{j}(C\cap D)
 = 
\epsilon^{\sigma}_{j, 1} (C)
=
\epsilon^{\sigma}_{j, 1}\circ \pi_{1}  (C, D)
$$ 
for
$
C\in\ch(\A_{1} )
$,
$
D\in\ch(\A_{2} )
$,
 and
$H_{j} \in \A_{1} $. 
\end{proof}

From now on, 
identify
$
\ch(\A_{1} )^{m} \times \ch(\A_{2} )^{m} 
$
and
$
\ch(\A_{1} \uplus \A_{2} )^{m} 
$
by the bijection $\alpha$ in Lemma \ref{Lemma3.1(1)}. 
Then Lemma \ref{Lemma3.1(2)} can be stated as
$$\epsilon_{j, p}^{\sigma}\circ \pi_{p} =
\epsilon_{j}^{\sigma}\,\,\,\,\,(p\in\{1, 2\}, \sigma\in B,
H_{j} \in \A_{p}).$$

\begin{proposition}
\label{Proposition12} 
There exists a natural bijection between
$AM(\A_{1} \uplus \A_{2} )$ and
$AM(\A_{1}) \times AM(\A_{2} )$.
\end{proposition}

\begin{proof}
Suppose that 
$\Phi$ is an admissible map 
for $\A_{1} \uplus \A_{2}$
and that 
a family of maps $\phhs$
$(H_{j} \in \A_{1} \uplus \A_{2}, \,\,\sigma\in B)$  
satisfies the conditions in Definition
\ref{admissiblemaps}.
Fix $p\in \{1, 2\}$
and $H_{j} \in \A_{p}$. 
Consider the following diagram:

\begin{diagram}
\ch(\A_{1})^{m} \times \ch(\A_{2})^{m}  & 
\rTo^{\Phi} & \ch(\A_{1}) \times \ch(\A_{2})\\
\dTo^{\pi_{p}} && \dTo^{\pi_{p}}\\
\ch(\A_{p})^{m} &\rDashto^{\Phi_{p} } &\ch(\A_{p})\\
\dTo^{\epsilon_{j,p}^{\sigma} } && \dTo^{\epsilon_{j,p}^{\sigma} }\\
B^{m} &
\rTo^{\varphi^{\sigma}_{j}} 
&
B\,\,.
\end{diagram}
By Lemma \ref{Lemma3.1(2)}, we have
$$
\epsilon_{j, p}^{\sigma}\circ \pi_{p}\circ \Phi 
=
\epsilon_{j}^{\sigma}\circ\Phi
=
\varphi^{\sigma}_{j}
\circ
\epsilon_{j}^{\sigma}
=
\varphi^{\sigma}_{j}
\circ
\epsilon_{j, p}^{\sigma}
\circ
\pi_{p} 
\,\,(p\in\{1, 2\}, \sigma\in B).$$ 
Assume
$p=1$ for simplicity.
Let $
C_{i} \in \ch(\A_{1} ),
D_{i} \in \ch(\A_{2} )
$ 
for 
$1\leq i\leq m$.
Then
\begin{align*} 
&~~~~\epsilon_{j, 1}^{\sigma}\circ \pi_{1}\circ \Phi 
(
C_{1} , C_{2} , \dots , C_{m}, 
D_{1} , D_{2} , \dots , D_{m})\\
&=
\varphi^{\sigma}_{j}
\circ
\epsilon_{j, 1}^{\sigma}
\circ
\pi_{1} 
(
C_{1} , C_{2} , \dots , C_{m}, 
D_{1} , D_{2} , \dots , D_{m})\\
&=
\varphi^{\sigma}_{j}
\circ
\epsilon_{j, 1}^{\sigma}
(
C_{1} , C_{2} , \dots , C_{m})
\end{align*} 
for each $H_{j} \in \A_{1} $.
Thus the chamber
$$\pi_{1}\circ \Phi 
(
C_{1} , C_{2} , \dots , C_{m}, 
D_{1} , D_{2} , \dots , D_{m})
\in
\ch(\A_{1})
$$
is independent of 
$D_{1} , D_{2} , \dots , D_{m}$.
Therefore we can express
  $$
\Phi_{1} (C_{1} , C_{2} , \dots , C_{m})
=
\pi_{1}\circ \Phi 
(
C_{1} , C_{2} , \dots , C_{m}, 
D_{1} , D_{2} , \dots , D_{m})
$$
for some map 
\[
\Phi_{1} : \ch(\A_{1} )^{m} \rightarrow 
 \ch(\A_{1} ).
\]
Then
$
\Phi_{1}
$ 
is an admissible map for $\A_{1}$
because the diagram above, including
$\Phi_{1} $,  is commutative
for each $H_{j}\in\A_{1}  $.
Simililarly we can define 
\[
\Phi_{2} : \ch(\A_{2} )^{m} \rightarrow 
 \ch(\A_{2} )
\]
so that
$
\Phi_{2}
$ 
is an admissible map for $\A_{2} $.
The construction so far gives a
natural map
\[
F : 
AM(\A_{1} \uplus \A_{2} )
\rightarrow
AM(\A_{1}) \times AM(\A_{2} ).
\]

Conversely suppose that 
$
\Phi_{p}$
is an admissible map for $\A_{p} $ 
and that 
a family of maps $\phhs$
$(H_{j} \in \A_{p}, \,\,\sigma\in B)$  
satisfies the conditions in Definition
\ref{admissiblemaps}.
Define
\[
\Phi:=\Phi_{1}\times \Phi_{2}
:
\ch(\A_{1})^{m} \times \ch(\A_{2})^{m}   
\longrightarrow \ch(\A_{1}) \times \ch(\A_{2}).
\]
Then
$
\Phi
$ 
is an admissible map for
$\A_{1}\uplus \A_{2}$
because the family of maps $\phhs$
$(H_{j} \in \A_{1}\uplus\A_{2},\,\,\sigma\in B)$  
satisfies the conditions in Definition
\ref{admissiblemaps}.
This construction gives a map
\[
G : 
AM(\A_{1}) \times AM(\A_{2} )
\rightarrow
AM(\A_{1} \uplus \A_{2} ).
\]
It is easy to check that $F$ and $G$ are 
inverses of each other.
\end{proof}

Now we have proved Theorem
\ref{maindecomp} by
applying Propositions \ref{decomp} and \ref{Proposition12}.

\section{Proof of Theorem \ref{main}} 
In this section we assume that
$\A
=
\{
H_{1} , H_{2} , \dots , H_{n} 
\}$ 
is a nonempty
real central {\it indecomposable} arrangement. 
We assume  $n\neq 2$ because any arrangement 
$\A = \{H_{1}, H_{2} \} $ is decomposable: 
\[
\A = \{H_{1} \}\uplus \{H_{2} \}.
\]

\begin{lemma}
\label{A1}
Let $m$ be a positive integer.
Suppose $\A$ is an arrangement with only one hyperplane $H_{1} $.
Let $H_{1} ^{+} $ and $H_{1} ^{-} $ be the two chambers.
Then

$(1)$ an arbitrary admissible map is given by
\[
\Phi(C_{1} , C_{2} , \dots , C_{m})
=
\begin{cases}
H_{1} ^{+}  \,\,\,\text{if~} C_{i} = H_{1} ^{+} \text{~for all~} i,\\
H_{1} ^{-}  \,\,\,\text{if~} C_{i} = H_{1} ^{-} \text{~for all~} i,\\
\text{either~} H_{1} ^{+} \text{~or~} H_{1} ^{-}  \,\,\,\text{~~~~otherwise},
\end{cases}
\]

$(2)$ the number of admissible maps is equal to 
$2^{2^{m}-2 }$,
and

$(3)$ every admissible map is projective if and only if
$m=1$. 
 
\end{lemma}

\begin{proof}
(1)
Note that the map
$\epsilon_{1}^{\sigma}
:
\ch(\A) \to B$ is a bijection.
So the commutativity condition 
in Definition \ref{admissiblemaps} 
can be ignored and
we simply consider a map $
{\Phi}
: \ch^{m} \rightarrow \ch$ satisfying
$
{\Phi}(H_{1} ^{\sigma}, H_{1} ^{\sigma}, \dots , H_{1} ^{\sigma}) = 
H_{1} ^{\sigma}
$ 
$
(\sigma \in B= \{+, -\}).
$ 

(2) We have two choices for 
each element of the set
\[
\ch^{m} \setminus \{(H_{1} ^{+}, H_{1} ^{+}, \dots , H_{1} ^{+}), 
(H_{1} ^{-}, H_{1} ^{-}, \dots , H_{1} ^{-})\}
\]
whose cardinality is equal to $2^{m} - 2$. 

(3) When $m=1$, 
by Proposition \ref{unanimity}, 
the only admissible map is the identity
map, which is projective. 
For $m\geq 2$, 
the number of admissible maps, which is equal to $2^{2^{m}-2 }$,
exceeds the number of projective ones, which is $m$. 
\end{proof}

Therefore we have proved Theorem \ref{main} (1).
Let us concentrate on
 Theorem \ref{main} (2).

\smallskip

Assume that $\A =
\{
H_{1}, H_{2} , \dots , H_{n}  
\}$ is indecomposable with $n=|\A|\geq 3$. 
Let $m$ be a positive integer.
We will show that every admissible map 
of $\A$ is projective.   
Suppose that
$\Phi$ is an admissible map
and that 
a family of maps $\phhs$
$(H_{j} \in \A, \,\,\sigma\in B)$  
satisfies the conditions in Definition
\ref{admissiblemaps}.

\begin{lemma}
\label{Lemma2.2B}
Assume
$1\leq j\leq n$
and 
$S\subseteq \{1, 2, \dots , m\}.$
Then
$
\varphi_{j}^{+} (S_{+}) 
=
-
\varphi_{j}^{-} (-S_{+}). 
$ 
In particular,  
$\varphi_{j}^{-} (-, -, \dots , -)
=
-. 
$ 
\end{lemma}

\begin{proof}
By Proposition \ref{surjective},  
we may choose
$
\mathcal C
\in\ch^{m} 
$ 
so that
$
\epsilon_{j}^{+} (\mathcal C) 
=
S_{+}.
$ 
Then
\begin{align*}
\varphi_{j}^{+}(S_{+} ) = +
&\Longleftrightarrow
\epsilon_{j}^{+}\circ \Phi(\mathcal C) 
= \varphi_{j}^{+}\circ\epsilon_{j}^{+} (\mathcal C) = +
\Longleftrightarrow
\Phi(\mathcal C) \subseteq H_{j}^{+}\\ 
&\Longleftrightarrow
-
=
\epsilon_{j}^{-}\circ \Phi(\mathcal C) 
= \varphi_{j}^{-}\circ\epsilon_{j}^{-} (\mathcal C)
=\varphi_{j}^{-}((S^{c})_{+})
=\varphi_{j}^{-}(-S_{+}).
\end{align*}
\end{proof}

Define 
$\delta_{\A}^{\sigma} : \ch(\A) \longrightarrow B^{n}$,
for $\sigma \in B = \{+, -\}$,  
by
\[
\delta_{\A}^{\sigma}(C) =
(
\epsilon_{1}^{\sigma}(C),
\epsilon_{2}^{\sigma}(C),
\dots,
\epsilon_{n}^{\sigma}(C)
).
\]
Then $\delta_{\A}^{\sigma}$ is injective.
We frequently suppress the subscript $\A$ in
$\delta_{\A}^{\sigma}$ 
when there is no fear of confusion.
Note that
$\delta^{+} (-C) = -\delta^{+}(C)
=
\delta^{-} (C)$, 
where $-C$ is the antipodal chamber of $C$.
 Thus $\delta^{-} = - \delta^{+} $.

\begin{lemma}
\label{Lemma2.2} 
Let 
$\B = 
\{
H_{1} , H_{2} , \dots , H_{\nu} 
\}
\subseteq
\A
$ be a circuit with $3\leq \nu\leq n$.
Then

$(1)$ $|\ch(\B)| = 2^{\nu} - 2$, and

$(2)$ there exists $\tau = (\tau_{1}, \tau_{2}, \dots , \tau_{\nu})\in B^{\nu}$
such that $$\im \delta_{\B}^{\sigma} = B^{\nu} \setminus 
\{\tau, -\tau\}.$$  
\end{lemma}

\begin{proof}
(1)
Since the intersection lattice 
\cite[Definition 2.1]{OrT}
$L(\B)$ of $\B$ 
is the same as that of the $\nu$-dimensional 
Boolean arrangement
(= the arrangement of the $\nu$ coordinate hyperplanes) 
in $\R^{\nu}$ up to the rank $\nu-1$,
the Poincar\'e polynomial
$\pi(\B, t)$ 
coincides with
the Poincar\'e polynomial of the 
$\nu$-dimensional Boolean arrangement
up to degree $\nu-1$.
The Poincar\'e polynomial of the 
$\nu$-dimensional Boolean arrangement
is equal to $(1+t)^{\nu}$
\cite[Example 2.49]{OrT}.  
Since $\deg \pi(\B, t) = r(\B)=
\nu-1$ and
$\pi(\B, -1) = 0$,
$\pi(\B, t) = (1+t)^{\nu} - t^{\nu} - t^{\nu-1} $. 
By \cite{Zas}
\cite[Theorem 2.68]{OrT}, one has
$|\ch(\B)|
=
\pi(\B, 1)
=
2^{\nu} - 2$.

(2)
By (1),
\[
\mid B^{\nu} \setminus \im \delta_{\B}^{+}  \mid
=
\mid B^{\nu} \mid - \mid\im \delta_{\B}^{+}  \mid
=
\mid B^{\nu} \mid - \mid \ch(\B)  \mid
=
2^{\nu} - (2^{\nu}-2)
=
2.
\]
Since $\delta_{\B}^{+} (-C)
=
-\delta_{\B}^{+} (C) $ for $C\in \ch(\B)$,
the set 
$\im \delta_{\B}^{+}$ is closed under the operation
$\tau\mapsto -\tau$. 
Thus
the set
$B^{\nu} \setminus \im \delta_{\B}^{+}$
is expressed as 
$\{\tau, -\tau\}$ for some $\tau\in B^{\nu}. $ 
\end{proof}

Define
\[
K_{j}^{\sigma} :=
\{
S\subseteq \{1, 2, \dots , m\}
\mid
\phis (S_{+})
=
+
\}
\,\,\,\,\,
(1\leq j\leq n,\,\, \sigma\in B= \{+, -\}).
\]

\begin{lemma}
\label{Lemma2.3}
Suppose that $\A$ is indecomposable and $n=|\A|\geq 3$. 
Then the maps
$\phis$ do not depend upon $j$
or $\sigma$.   
\end{lemma}

\begin{proof}
Choose a circuit
$\B\subseteq \A$.  We may assume that
$\B = 
\{
H_{1} , H_{2} , \dots , H_{\nu} 
\}
$ and $3\leq \nu\leq n$.
By Lemma \ref{Lemma2.2},
there exists 
$
\tau 
=
(
\tau_{1},
\tau_{2},
\dots,
\tau_{\nu}
)
\in B^{\nu} $ such that
\[
B^{\nu} 
=
(\im \delta_{\B}^{+})
\cup
\{
\tau,
-\tau
\}
\,\,\,
\text{(disjoint)}.
\]
Let 
$1\leq p\leq \nu$,
$1\leq q\leq \nu$,
$p\neq q$.
Since
neither of
$
(
\tau_{1},
\dots,
-\tau_{q},
\dots,
\tau_{\nu}
)
$
nor
$
(
\tau_{1},
\dots,
-\tau_{p},
\dots,
\tau_{\nu}
)
$
lies in 
$\{
\tau,
-\tau
\}
$,
they both lie in $\im\delta_{\B}^{+} $.
Choose
$C, C'\in\ch(\B)$ 
such that
$$
\delta_{\B}^{+}(C) =
(
\tau_{1},
\dots,
-\tau_{p},
\dots,
\tau_{\nu}
)
,
\,\,\,
\delta_{\B}^{+}(C') =
(
\tau_{1},
\dots,
-\tau_{q},
\dots,
\tau_{\nu}
).
$$   
Choose 
$\hat{C}\in\ch(\A) $ 
and
$\hat{C'}\in\ch(\A) $ 
so that
$\hat{C}\subseteq C$
and
$\hat{C'}\subseteq C'$. 
Let $S\subseteq
\{
1, 2, \dots, m
\}
$. 
Define
$
\mathcal C
=
(
C_{1} , C_{2} , \dots , C_{m} 
)\in\ch(\A)^{m} 
$  by
\[
C_{i} 
=
\begin{cases}
\hat{C'} \,\,\text{if~} i\in S,\\
\hat{C} \,\,\text{if~} i\not\in S.
\end{cases}
\]
Then
\[
\epsilon_{p}^{\tau_{p} } (\mathcal C) 
=
\epsilon_{q}^{-\tau_{q} } (\mathcal C) 
=
S_{+} ,
\,\,\,
\epsilon_{r}^{\tau_{r} } (\mathcal C) 
=
(+, +, \dots , +)
\,\,\,
(1\leq r \leq \nu, r \not\in\{p, q\}).
\]
Suppose 
$S\in K_{p}^{\tau_{p}}$, 
i. e.,
$
\varphi_{p}^{\tau_{p} } 
(S_{+} ) = +
$.
Then
\[
\epsilon_{p}^{\tau_{p} } 
\circ
\Phi
(\mathcal C) 
=
\varphi_{p}^{\tau_{p} }
\circ
\epsilon_{p}^{\tau_{p} } (\mathcal C) 
=
\varphi_{p}^{\tau_{p} } 
(S_{+} ) = +.
\]
This implies that
$\Phi(\mathcal C) \subseteq 
H_{p}^{\tau_{p} }.$ 
Similarly we have
$
\Phi(\mathcal C) \subseteq 
H_{r}^{\tau_{r} } 
$ 
when $
1\leq r \leq \nu, r \not\in\{p, q\}
$, 
because
$
\varphi_{r}^{\tau_{r} }
\circ
\epsilon_{r}^{\tau_{r} } (\mathcal C) 
=
\varphi_{r}^{\tau_{r} }
(+, +, \dots, +) = +.
$ 
Note that
\[
\bigcap_{j=1}^{\nu}  H_{j}^{\tau_{j} } = \emptyset 
\]
because $\tau \not\in \im \delta_{\B}^{+} $. 
Therefore
\[
\Phi(\mathcal C) \subseteq
\bigcap_{j\neq q}^{\nu}  H_{j}^{\tau_{j} }
\subseteq
H_{q}^{-\tau_{q} }.  
\]
Thus
\[
\varphi_{q}^{-\tau_{q} }
(S_{+})
=
\varphi_{q}^{-\tau_{q} }
\circ
\epsilon_{q}^{-\tau_{q} } (\mathcal C) 
=
\epsilon_{q}^{-\tau_{q} } 
\circ
\Phi
(\mathcal C) 
=
+,
\]
which implies
$S\in 
K_{q}^{-\tau_{q} }
$. 
Therefore
$
K_{p}^{\tau_{p} } 
\subseteq
K_{q}^{-\tau_{q} }. 
$ 

Similarly one can show
$
K_{p}^{\tau_{p} } 
\supseteq
K_{q}^{-\tau_{q} }, 
$ and thus
$
K_{p}^{\tau_{p} } 
=
K_{q}^{-\tau_{q} }
$ if $p\neq q$. 
Since $\nu
\geq 3$,
 we can conclude that
$K_{j}^{\sigma}$ does not depend  
upon
$j$, 
$1\leq j\leq \nu$, or $\sigma\in B$. 
So 
$\varphi _{j}^{\sigma}$ does not depend  
upon
$j$,  
$1\leq j\leq \nu$, or $\sigma\in B$. 
Apply Lemma \ref{lemma2.1}, and we know
$\varphi _{j}^{\sigma}$ does not depend  
upon
$j$, 
$1\leq j\leq n$, or $\sigma\in B$. 
\end{proof}

Because of Lemma \ref{Lemma2.3},
write $\varphi=\phhs$ 
for $j$,
$1\leq j\leq n$,
 and $\sigma\in B$.
Let
$$K 
=
\{
S\subseteq \{1,2,\dots,m\}
\mid
\varphi(S_{+} ) = +
\}.$$

\begin{lemma}
\label{Lemma2.4}

$(1)$
$\{1, \dots , m\}\in K$,
$(2)$
$S\in K$ if and only if
$S^{c} \not\in K$,
$(3)$
 $S_{1} \cap S_{2} \in K$ if
 $S_{1} \in K$ and
$S_{2} \in K$.
\end{lemma}

\begin{proof}
(1) is obvious because
$\varphi (+, +, \dots , +)=+$. 

(2) By Lemma \ref{Lemma2.2B} 
\begin{align*}
S\in K = K_{1}^{+}  
&\Longleftrightarrow
\varphi_{1}^{+}(S_{+} ) = +
\Longleftrightarrow
\varphi_{1}^{-}((S^{c})_{+}) 
=
\varphi_{1}^{-}(-S_{+}) 
=
-
\varphi_{1}^{+}(S_{+} ) = 
 -\\
&\Longleftrightarrow
\varphi_{1}^{-}((S^{c})_{+}) 
=-
\Longleftrightarrow
S^{c} \not\in K_{1}^{-}=K.
\end{align*}

(3) Choose a circuit $\B \subseteq \A$. 
We may assume 
$\B =
\{
H_{1} , H_{2} , \dots , H_{\nu} 
\}$
with $3\leq \nu \leq n$.
By Lemma \ref{Lemma2.2},
there exists 
$
\tau 
=
(
\tau_{1},
\tau_{2},
\dots,
\tau_{\nu}
)
\in B^{\nu} $ such that
\[
B^{\nu} 
=
(\im \delta^{+})
\cup
\{
\tau,
-\tau
\}
\,\,\,
\text{(disjoint)}.
\]
There exist four chambers
$C, C', C'', C'''\in \ch(\B)$ such that
\begin{align*} 
\delta_{\B}^{+}(C) 
&=
(
\tau_{1},
\tau_{2},
-\tau_{3},
\tau_{4},
\dots,
\tau_{\nu}
)
,
\,\,\,
\delta_{\B}^{+}(C') =
(
\tau_{1},
-\tau_{2},
\tau_{3},
\tau_{4},
\dots,
\tau_{\nu}
),\\
\delta_{\B}^{+}(C'') 
&=
(
-\tau_{1},
\tau_{2},
\tau_{3},
\tau_{4},
\dots,
\tau_{\nu}
)
,
\,\,\,
\delta_{\B}^{+}(C''') =
(
-\tau_{1},
-\tau_{2},
\tau_{3},
\tau_{4},
\dots,
\tau_{\nu}
).
\end{align*} 
Choose four chambers 
$
\hat{C}    ,
\hat{C'},
\hat{C''},
\hat{C'''}
\in\ch(\A)  
$ 
such that
$$
\hat{C} \subseteq C,\,\,\,
\hat{C'} \subseteq C',\,\,\,
\hat{C''} \subseteq C'',\,\,\,
\hat{C'''} \subseteq C'''.
$$
Assume that $S_{1}, S_{2} \in K
$. 
Define
$
\mathcal C
=
(
C_{1} , C_{2} , \dots , C_{m} 
)\in\ch(\A)^{m} 
$  by
\[
C_{i} 
=
\begin{cases}
\hat{C}  \,\,\text{if~} i\in S_{1} \cap S_{2},\\
\hat{C'}  \,\,\text{if~} i\in S_{1} \setminus S_{2},\\
\hat{C''}  \,\,\text{if~} i\in S_{2} \setminus S_{1},\\
\hat{C'''}  \,\,\text{if~} i\not\in S_{1} \cup S_{2}.
\end{cases}
\]
Then
\begin{align*} 
\epsilon_{1}^{\tau_{1} } (\mathcal C) 
&=
(S_{1})_{+},
\,\,\,
\epsilon_{2}^{\tau_{2} } (\mathcal C) 
=
(S_{2})_{+},
\,\,\,
\epsilon_{3}^{-\tau_{3} } (\mathcal C) 
=
(S_{1}\cap S_{2})_{+},\\
\epsilon_{j}^{\tau_{j} } (\mathcal C) 
&=
(+, +, \dots, +)
\,\,\,
(4\leq j\leq \nu).
\end{align*} 
Thus we have
\begin{align*} 
\epsilon_{1}^{\tau_{1} } \circ \Phi(\mathcal C) 
&=
\varphi \circ
\epsilon_{1}^{\tau_{1} }(\mathcal C) 
=
\varphi ((S_{1})_{+}) 
=+,\\
\epsilon_{2}^{\tau_{2} } \circ \Phi(\mathcal C) 
&=
\varphi \circ
\epsilon_{2}^{\tau_{2} }(\mathcal C) 
=
\varphi((S_{2})_{+}) 
=+,\\
\epsilon_{j}^{\tau_{j} } \circ \Phi(\mathcal C) 
&=
\varphi \circ
\epsilon_{j}^{\tau_{j} }(\mathcal C) 
=
\varphi(+, +, \dots ,+) 
=+
\,\,\,
(4\leq j\leq \nu),
\end{align*} 
which implies 
\[
\Phi(\mathcal C)
\subseteq
H_{1}^{\tau_{1}} 
\cap
H_{2}^{\tau_{2}} 
\cap
H_{4}^{\tau_{4}} 
\cap
\dots
\cap
H_{\nu}^{\tau_{\nu}} 
\subseteq
 H_{3}^{-\tau_{3}}.  
\]
Therefore 
\[
\varphi 
((S_{1} \cap S_{2})_{+})
=
\varphi 
\circ
\epsilon_{3}^{-\tau_{3} }  
(\mathcal C)
=
\epsilon_{3}^{-\tau_{3} }  
\circ
\Phi
(\mathcal C)
=
+
\]
and 
$S_{1} \cap S_{2} \in K$. 
\end{proof}

Now we are ready to prove the following statement,
which is Theorem \ref{main} (2).

\begin{quotation} 
\noindent
{\it
Let $\A$ be a real central indecomposable 
arrangement with 
$|\A| \geq 3$. Then
every admissible map is projective.
}
\end{quotation} 

\begin{proof}
Define
$S_{0} = \bigcap_{S\in K} S$.  By Lemma \ref{Lemma2.4} (3),
$S_{0} \in K$.  By Lemma   \ref{Lemma2.4} (1) and (2),
we have
$\emptyset\not\in K$.  Thus $S_{0} \neq \emptyset$.
Let $h\in S_{0} $.
Since $S_{0} \setminus \{h\}\not\in K$,
$(
\{1, 2, \dots , m\}
\setminus S_{0})\cup\{h\} \in K$ by
Lemma   \ref{Lemma2.4} (2).
By Lemma \ref{Lemma2.4} (3),
\[
\{h\} = ((
\{1, 2, \dots , m\}
\setminus S_{0})\cup\{h\}) \cap S_{0} \in K.
\]
Thus $S_{0} = \{h\}$. 
Note that, by 
Lemma   \ref{Lemma2.4} (2),
$$
S\in K
\Rightarrow
h\in S
\Leftrightarrow
h\not\in S^{c}
\Rightarrow
S^{c} \not\in K
\Leftrightarrow
S\in K.
$$ 
Therefore,
$
S\in K
$
if and only if
$h\in S
$:
\[
K = \{S \subseteq 
\{1, 2, \dots , m\}
\mid h\in S\}.
\]
This implies that $\varphi $ is equal to the
projection to the $h$-th component.
Let 
$\mathcal C\in \ch^{m} $.
Then
\begin{align*}
\epsilon_{j}^{\sigma}\circ \Phi(\mathcal C)
=
\varphi \circ \epsilon_{j}^{\sigma} (\mathcal C)
=
\varphi (
\epsilon_{j}^{\sigma} (C_{1}),
\epsilon_{j}^{\sigma} (C_{2}),
\dots,
\epsilon_{j}^{\sigma} (C_{m}))
=
\epsilon_{j}^{\sigma} (C_{h}).
\end{align*}
  Since $\Phi(\mathcal C)$ and 
$C_{h} $ lie on the same side of 
every hyperplane $H_{j} \in \A$,
 $\Phi(\mathcal C)
=
C_{h}.$  Therefore
 $\Phi$ is the projection to the $h$-th component. 
\end{proof}

Decompose a nonempty
real central arrangement $\A$ into
nonempty
indecomposable arrangements as
$$
\label{indecomposabledecomposition}
\A = 
\A_{1} \uplus \A_{2} \uplus \dots\uplus\A_{a} 
\uplus
\B_{1} \uplus \B_{2} \uplus \dots\uplus\B_{b},
\eqno{(**)}
$$
where 
$|\A_{p}| = 1 \,\, (1\leq p\leq a)$ 
and
$|\B_{q}| \geq 3 \,\, (1\leq q\leq b)$.
Then,
by Lemma \ref{A1}, Theorems \ref{maindecomp} and \ref{main},  the number of admissible maps for $\A$ is equal to
\[
\left(2^{2^{m}-2 }\right)^{a} m^{b}.  
\]
This proves
Corollary \ref{corollary1}.

Next we will prove
Corollary \ref{corollary2}:
If $m=1$, then,
by Proposition \ref{unanimity},
  the only admissible map is 
the identity map $\ch \to \ch$, which is projective.  
Assume $m\geq 2$. Then,
by Lemma \ref{A1}, Theorems \ref{maindecomp} and \ref{main},  
every admissible map is projective
if and only if 
$a=0$ and $b=1$ 
in the decomposition
$(**)$
 above.

\bigskip

{\small {\bf Acknowledgement.} 
The author would like to express his gratitude to
Professors H. Kamiya and A. Takemura, who introduced
him to Arrow's impossibility theorem and gave
him helpful comments for  earlier versions of this paper, 
and to Dr. T. Abe
with whom he had stimulating conversations.}

\end{document}